\documentstyle[12pt]{article}
\title{ ZERO ASSIGNMENT IN MULTIVARIABLE SYSTEM \\
           USING POLE ASSIGNMENT METHOD }
\author{\bf Ye.M.\,Smagina\\[2ex]
        \it Department of Chemical Engineering, Technion, \\
        \it Haifa, Israel 32000\\
        \it email: cermssy@tx.technion.ac.il%
       }
\date{}
\newcommand{\be}{\begin{equation}}
\newcommand{\ee}{\end{equation}}

\newcommand{\ba}{\begin{array}}
\newcommand{\ea}{\end{array}}

\newtheorem{ass}{Assertion}
\newtheorem{rem}{Remark}
\newtheorem{pro}{Problem}
\newcommand{\rf}[1]{{\rm(\ref{#1})}}
\newtheorem{th}{Theorem}
\newcommand{\qed}{\mbox{Q.E.D.}}

\begin{document}

\maketitle
\pagestyle{empty}

\begin{abstract}
In the paper we consider the invariant zero assignment problem in
a linear multivariable system with several inputs/outputs by
constructing a system output matrix. The problem is reduced to the pole
assignment problem by a state feedback (modal control) in a descriptor
system or a regular one. It is shown that the zero assignment and pole
assignment are mathematically equivalent problems.
\end{abstract}

\hyphenation{ }

\section{\small{INTRODUCTION}}
For the first time the problem of the arbitrary zero assignment by
choosing the output matrix was studied in (Rosenbrock, 1970, Theorem
4.1) and  it is formulated briefly as follows:
let an output matrix be arbitrarily chosen for a system.
It is necessary to find this matrix so that the obtained system is observable
and its transfer function matrix  has the assigned numerators $\epsilon_i(s)$ of
diagonal elements of Smith-McMillan form.

Rosenbrock showed that such the output matrix can always be chosen
if the degrees of $\epsilon_i(s)$ satisfy some conditions. But these
conditions are rather complicate for the designing. From this point of view
we will consider the zero assignment method which assigns the zero
polynomial $\psi(s) = \epsilon_1(s)\epsilon_2(s) \cdots \epsilon_r(s)$
entirely. Such a statement simplifies considerably Rosenbrock's conditions.

The problem of the zero polynomial assignment (zero assignment) has been
studied in the works of the author. In the first works (Smagina,1984;
Smagina,1986) the iterative zero assignment method was  proposed,
in the following works (Smagina,1985; Smagina,1991) the analytic zero
assignment method was studied and it was discovered (Smagina, 1985)
that the zero assignment problem is mathematically equivalent to the pole
assignment problem. Later on the above result was formally proved using the
special Yokoyama's canonical form (Smagina, 1996). A similar result has been
obtained in (Syrmas and Lewis,1993) by using Hessenberg's canonical form.

In the present paper we prove the equivalence between the zero and pole
assignment problems without using the complex canonical transformations.
It is shown that, in general, the zero assignment is equivalent to the pole
assignment in a descriptor system. When there is no restriction on an
assigned zero polynomial degree then the zero assignment is
equivalent to the pole assignment in a reduced regular dynamical system.

\section{\small{PROBLEM STATEMENT}}
We consider a linear multivariable dynamical system in the state-space

\be
   \dot{x} = Ax+Bu,
\label{eq1}
\ee
with the output
\be
  y = Hx
\label{eq2}
\ee
where $x$ is a state $n$-vector, $u$ is an input $r$-vector, $y$ is an output
$l$-vector $(r,l \le n)$; $A,B,H$ are real constant matrices of the
appropriate sizes.  It is assumed that $rankB=r$, $rankH =l$.

 It is known (MacFarlane and Karcanias, 1976) that finite invariant zeros of
system (1), (2) are defined as the set of complex $s$ which satisfies the
rank equality
\be
rankP(s) = rank\left[ \begin{array}{cc} sI-A & -B \\ H  &
O \end{array} \right] < n + min(r,l)
\label {eq3}
\ee
If $r=l$ then the finite invariant zeros are defined as zeros of the
polynomial
\be
\psi(s) = detP(s)
\label{eq4}
\ee
Since there exist no invariant zeros in "almost all" system (1), (2) with
$r\ne l$ (Davison and Wang, 1974) then we will consider the zero assignment
problem for the "square" system \rf{eq1},\rf{eq2} with $r = l$. \footnote
{ If $r \neq l$ then it would be desirable to solve the zero assignment
problem by "squaring down" system inputs or outputs.}
Therefore, the zero assignment is formulated as follows.
\begin{pro}\label{P1}
For the pair
$(A,B)$ find an output matrix $H$ which ensures the coincidence of the
polynomial $\psi(s)$ with a preassigned polynomial
\be \psi_a(s) = \beta_o
+\beta_1s +\cdots+ \beta_{\mu}s^{\mu}
\label{eq5}
\ee
where $\mu \le n-r$.
\footnote{The maximum number of zeros in a $n$  order system with $r$
inputs/outputs is $n-r$ (Davison and Wang, 1974).} \end{pro}

From the practical motivation it is necessary to impose the following condition
on the matrix $H$:

\be
i) \rm {\;\;the \;pair \;}(A,H)\; \rm{is \; observable}.
\label{eq6}
\ee

\begin{rem}\label{R1}
It has been shown in the previous works of the author (Smagina, 1984;
Smagina, 1986) that condition \rf{eq6} is fulfilled if any zero
$\bar{s}_i, i=1,\ldots,\mu$ of the preassigned polynomial $\psi_a(s)$ \rf{eq5}
does not coincide with any eigenvalue $\lambda_j(A), j=1,\ldots,n$ of the
matrix $A$
\be \bar{s}_i \neq \lambda_j(A), \;\;i=1,\ldots,\mu, j=1,\ldots,n.
\label{eq7}
\ee
\end{rem}

\section{\small{REDUCTION TO POLE ASSIGNMENT IN DESCRIPTOR SYSTEM}}

 It follows from the condition $rank B = r$ that $B$ has $r$ linearly
independent rows.  At first for simplicity we consider the matrix $B$ of the
following structure
\be
    B=\left[ \begin{array}{r} B_1 \\B_2
       \end{array} \right]
\label{eq8}
\ee
where $B_1$ and $B_2$ are $(n-r)\times r$ and $r\times r$ blocks respectively,
moreover,
$$
 detB_2 \neq 0.
$$
We introduce the $n\times n$-nonsingular transformation matrix
\be
    N=\left[ \begin{array}{cc}
              I_{n-r} & -B_1B_2^{-1} \\
                O  &  I_r
       \end{array}           \right]
\label{eq9}
\ee
which reduces the matrix $B$ to the form
$$
    NB=\left[ \begin{array}{cc}
              I_{n-r} & -B_1B_2^{-1} \\
                O  &  I_r
       \end{array} \right] \left[ \begin{array}{r} B_1 \\B_2
        \end{array} \right] \;=\; \left[ \begin{array}{r} O\\B_2
        \end{array} \right]
$$
and we transform the system matrix $P(s)$ into $\bar{P}(s)$ by means of the
strict system equivalence operations (Pugh et al, 1987) which preserve
the finite and infinite zeros of the elementary divisors of a regular matrix
pencil:

\be \bar{P}(s) = \left[ \begin{array}{cc} N & O \\ O & I_l\end{array} \right]
        P(s)\left[ \begin{array}{cc} N^{-1} & O \\ O & I_r\end{array} \right]
\;\;=\;\; \left[ \begin{array}{cc} sI-NAN^{-1} & -NB \\ HN^{-1} &  O
          \end{array} \right]
\label{eq10}
\ee
 We decompose the matrices $A$ and $H$ as
$$
     A \;=\;\left[ \begin{array}{cc} A_{11} & A_{12} \\ A_{21}& A_{22}
            \end{array} \right], \qquad
             H \;=\; \left[ \begin{array}{cc} H_1 & H_2 \end{array} \right]
$$
with the blocks $A_{11}, A_{12}, A_{21}, A_{22}, H_1, H_2$ having sizes
$(n-r)\times (n-r)$, $(n-r)\times r$, $r\times (n-r)$, $r\times r$, $l\times
 (n-r)$, $l\times r$ respectively and calculate

$$
 NAN^{-1} = \left[ \begin{array}{cc} \bar{A}_{11} & \bar{A}_{12} \\
           \bar{A}_{21}& \bar{A}_{22} \end{array} \right ],\;\;
   HN^{-1}  = \left [ \begin{array}{cc} \bar{H}_1 & \bar{H}_2 \end{array}
           \right ]
$$

where

$$
\bar{A}_{11}= A_{11} -B_1B_2^{-1}A_{21},\;\;
$$
$$
\bar{A}_{12}=(A_{11}-B_1B_2^{-1}A_{21})B_1B_2^{-1}+A_{12}-B_1B_2^{-1}A_{22},
$$
$$
  \bar{A}_{21} = A_{21}, \;\;  \bar{A}_{22} = A_{21}B_1B_2^{-1} +A_{22},
$$
\be
     \bar{H}_1 = H_1,\;\;   \bar{H}_2 = H_2B_1B_2^{-1} +H_2
\label{eq11}
\ee

Using (11), we represent $\bar{P}(s)$ as follows

\be
\bar{P}(s) = \left[ \begin{array}{ccc} sI_{n-r}-\bar{A}_{11} & -\bar{A}_{12} &
O\\ - \bar{A}_{21}& sI_r - \bar{A}_{22} & -B_2 \\ \bar{H}_1 & \bar{H}_2 & O
 \end{array} \right]
\label{eq12}
\ee

The finite and infinite zeros (Verghese et al.,1981) of $P(s)$ and
$\bar{P}(s)$ coincide because the strict system equivalence operations
preserve the finite and infinite elementary divisors of the regular pencil.
Hence, the finite zeros of the system \rf{eq1},\rf{eq2} with $r = l$
coincide with zeros of the polynomial $det\bar{P}(s)$:
\be
\psi(s) = detP(s) = det\bar{P}(s).
\label{eq13}
\ee

Using the block structure of $\bar{P}(s)$, we calculate
\be
det\bar{P}(s) = (-1)^q detB_2det \left[ \begin{array}{cc}
   sI_{n-r}-\bar{A}_{11} & -\bar{A}_{12}\\ \bar{H}_1 & \bar{H}_2
 \end{array} \right]
\label{eq14}
\ee
where an integer $q > 0$.

As a result {\it Problem~\ref{P1}} is reformulated as follows.
\begin{pro}\label{P2}
Find matrices $\bar{H}_1$ and  $\bar{H}_2$ such that the following equality
could take place
\be
det\left[ \begin{array}{cc} sI_{n-r}-\bar{A}_{11} & -\bar{A}_{12}\\
       \bar{H}_1 & \bar{H}_2  \end{array} \right] \;=\;\psi_a(s)
\label{eq15}
\ee
\end{pro}

Condition \rf{eq6} will be guaranteed by \rf{eq7}.

If  $\bar{H}_1$ and $\bar{H}_2$ are found, the original output matrix
$ H = \left[ \begin{array}{cc} H_1 & H_2 \end{array} \right]$ of system
\rf{eq1},\rf{eq2}  is defined as follows
$$
H = \left[ \begin{array}{cc} H_1 & H_2 \end{array} \right] =
    \left[ \begin{array}{cc} \bar{H}_1 & \bar{H}_2 \end{array} \right]N =
    \left[ \begin{array}{cc} \bar{H}_1 & -\bar{H}_1B_1B_2^{-1} + \bar{H}_2
     \end{array} \right].
$$
Therefore,
\be
   H_1 = \bar{H}_1,\;\; H_2 = \bar{H}_2 - \bar{H}_1B_1B_2^{-1}.
\label{eq16}
\ee

Now we will show that {\it Problem~\ref{P2}} is equivalent to a pole assignment
  problem in a descriptor system formed from the blocks of the matrix in
\rf{eq15}.

Denote
$$
 E = \left[ \begin{array}{cc} I_{n-r} & O \\ O & O \end{array} \right],\;
 F = \left[\begin{array}{cc} \bar{A}_{11} & \bar{A}_{12} \\O& O
\end{array}\right],\; G = \left[\begin{array}{c} O \\ -I_r \end{array}\right],\;
$$
\be
             K \;=\; \left[ \begin{array}{cc} \bar{H}_1 & \bar{H}_2
\end{array} \right]
\label{eq17}
\ee
\begin{th}\label{Th1}
The zero assignment problem in system \rf{eq1},\rf{eq2} with $rankB =r$ is
equivalent to the pole assignment problem in the following descriptor system
\be
            E\dot{w} = Fw +Gv
\label{eq18}
\ee
by using the feedback proportional state regulator
\be
   v = Kw
\label{eq19}
\ee
with the $r\times n$ matrix $K$.
\end{th}

Therefore, the matrix $K$ obtained will guarantee the condition
\be
det(sE - F - GK) = \psi_a(s)
\label{eq20}
\ee
In \rf{eq18} $w$ is a state $n$ vector, $v$ is an input $r$ vector.

\it{Proof}. \rm The proof  follows directly from \rf{eq4}, \rf{eq13},
\rf{eq15} and notations \rf{eq17} i.e.
 $$
         detP(s) = det\bar{P} = det\left[ \begin{array}{cc}
  sI_{n-r}-\bar{A}_{11} & -\bar{A}_{12}\\ \bar{H}_1 & \bar{H}_2  \end{array}
         \right] =  det(sI - F + GK)\;\;\qed
$$

In the remaining part of this section we consider the matrix $B$ having an arbitrary
structure. By applying the state transformation $\bar{x} = Mx$ \footnote{ $M$
is a permutation $n\times n$ matrix which has a single unity element in each
row(column) and zeros otherwise.} which rearranges the rows of $B$ in
\rf{eq1} we bring the matrix $B$ to the form
\be
 MB = \left [ \begin{array}{c} \tilde{B}_1\\ \tilde{B}_2 \end{array} \right]
\label{eq21}
\ee
with $rank\tilde{B}_2 =r$. Then we consider the transformation matrix which
has structure \rf{eq9} and brings the matrix $P(s)$ to $\bar{P}(s)$
by strict system equivalence transformations \rf{eq10} with $N$
replaced by
\be \bar{N} = NM.
\label{eq22}
\ee
Further, we can use the above approach with $\bar{H} = [\bar{H}_1,\bar{H}_2]$
$ = H\bar{N}^{-1}$. Then from the relation $ H = \bar{H}NM $ we can
calculate the matrix $H$ as follows
\be
H = [\bar{H}_1,\bar{H}_2] \left[ \begin{array}{cc}
I_{n-r} & -\tilde{B}_1\tilde{B}_2^{-1} \\ O  &  I_r \end{array}
\right] M = \left [ \bar{H}_1, \;\bar{H}_2 - \bar{H}_1\tilde{B}_1
\tilde{B}_2^{-1}\right ] M
\label{eq23}
\ee
Therefore, we have shown that the zero assignment problem in system
\rf{eq1}, \rf{eq2} is reduced to the pole assignment problem in the
descriptor system  \rf{eq18} of order $n$ with $r$ inputs. Recently a number
of authors (Armantano, 1984; Chu, 1988; Blanchini,1990; etc.) has considered
the pole assignment problem by the proportional state feedback in a descriptor
system of the general structure. The pole assignment solvability conditions
are simplified in system \rf{eq18} in virtue of the special structure of
the matrices $E$, $F$, $G$.
Further we will analyse the solvability conditions in detail and reveal their
connection with the structure of original system \rf{eq1}, \rf{eq2}.

\section{\small{ANALYSIS OF POLE ASSIGNMENT SOLVABILITY CONDITIONS
IN DESCRIPTOR SYSTEM}}

The pole assignment in descriptor system \rf{eq18} may be formulated as
follows:

Find a state feedback control law \rf{eq19} so that

i) the closed loop pencil $F_c(s) = sE -(F + GK)$ is regular: $detF_c(s)
   \neq 0$ (Gantmacher, 1959);

ii) the regular pencil $F_c(s)$ has $\mu \leq n-r$ preassigned (finite)
    eigenvalues.

We use the definition and the theorem of the regularizability  from
(Ozcaldizan and Lewis, 1990) in order to analyse the point i).
System \rf{eq18} is said to be regularizable by a proportional state
feedback if the pencil $(sE -(F+GK)$ is regular for some matrix $K$.
System \rf{eq18} is always regularizable if and only if
\be
dim(E{\cal V}^* + Im G) = dim {\cal V}^*
\label{eq24}
\ee
where ${\cal V}^*$ is a supremal element of the $(F,E,G)$-invariant subspace
$\cal V$ : $F{\cal V} \in E{\cal V} + Im G$ where $Im G$ is the
image of $G$.  It is clear that relation \rf{eq24} holds for system
\rf{eq18} in virtue of the special structure of $E$, $F$, $G$ \rf{eq17}.
Therefore, the condition i) is always satisfied for system \rf{eq18}.

Then we analyse the point ii). As it is known the regular pencil $F_c(s)$ has
$q=rankE$ eigenvalues (finite and infinite). It follows from the works
(Cobb,1981; Verghese et al,1981; Yip and Sincovec,1981) that only
controllable finite and infinite eigenvalues may be arbitrarily assigned.
The finite eigenvalues of $F_c(s)$ are controllable if and only if
$$ rank(sE -(F+GK), G) =n. $$
Since $ rank(sE -(F+GK), G) = rank(sE - F, G)$ then we have the
following controllability condition of the finite eigenvalues of $F_c(s)$
\be
     rank(sE -F, G) =n.
\label{eq25}
\ee
The infinite eigenvalues of $F_c(s)$ are controllable and may be shifted to
preassigned points of the complex plan if the following condition takes place
(Armantano, 1984; Blanchini,1990)
\be
     rank(E,G) = n. \footnote{ Condition: $rank[E, F{\cal S}_o, G] =n$
(Armantano, 1984) takes place if \rf{eq26} is satisfied
where ${\cal S}_o = kerE$.}
\label{eq26}
\ee
Therefore, condition \rf{eq25} and  \rf{eq26}
guarantee that a proportional feedback $r\times n$ matrix $K$ exists so that
the pole assignment problem has a solution.

It is evident that condition \rf{eq26} is always fulfilled if $rankB =r$ by
virtue of the special structure of $E$ and $G$.
The analysis of \rf{eq25} gives
\begin{ass}\label{A1}
The finite eigenvalues of the pencil $F_c(s)$ are controllable if and only if
the pair $(A,B)$ is controllable,
\end{ass}

\it{Proof}. \rm The proof follows from the rank equalities:
$$
n = rank\left [ sI - A, \; B \right ] = rank\left ( \left [ \begin{array}{cc}
\bar{N} & O \\ O & I  \end{array} \right ] \left[ sI - A, \; B \right ]
\bar{N}^{-1} \right ) = $$
$$  = rank \left [ \begin{array}{ccc} sI-\tilde{A}_{11} &
  -\bar{A}_{12} & O \\ -\bar{A}_{21} & sI-\bar{A}_{22} & \tilde{B}_2
   \end{array} \right ]  =
    rank \left [ sI - \bar{A}_{11}, \; - \bar{A}_{12} \right ]  + r = $$
$$    rank \left [ \begin{array}{ccc} sI-\bar{A}_{11} &
 -\bar{A}_{12} & O \\ O & O & I_r  \end{array} \right]  =
  rank\left [ sE - F, \; G \right ]\qed
$$
Therefore, the solvability conditions of the pole assignment problem in
descriptor system \rf{eq18}  are the controllability of the pair $(A,B)$
and the rank fullness of $B$. Simultaneously these conditions are the
solvability conditions of the zero assignment problem in system \rf{eq1},
\rf{eq2}.  Taking into account {\it Remark~\ref{R1}}  we have

\begin{th}\label{Th2} If the pair $(A,B)$ is controllable, $rankB =r$ and zeros
$z_i$, $i=1,\cdots, \mu$ of a preassigned polynomial $\psi_a(s)$ don't coincide
with eigenvalues of $A$, then there always exists the $r\times n$ matrix $H$
which ensures that system \rf{eq1}, \rf{eq2} is observable and has the
preassigned zeros $z_i$, $i=1,\cdots,\mu$, $\mu \leq n-r$.  \end{th}

{\it Theorem~\ref{Th2}} generalizes the appropriate results of the
previous works of the author (Smagina, 1985; Smagina, 1996).

\section{\small{REDUCTION TO POLE ASSIGNMENT IN REGULAR STATE-SPACE SYSTEM}}

It is known (Kouvaritakis and MacFarlane,1976) that the maximal number of
finite invariant zeros in a $n$-order system with $r$ inputs and outputs is
equal  to $n-r$. Using the method of sect.4 we can arbitrarily assign any number
$\mu \leq n-r$ of finite zeros (remaining $n-r-\mu$ zeros are situated in
$\infty$). If there is no the restriction on a number of finite zeros,
namely, we can assign $\mu = n-r$ then {\it Problem~\ref{P1}} can be
reduced to the pole assignment in a regular state-space system
(regular pole-assignment problem). Indeed, the output $r\times n$ matrix $H$
has $rn$ elements, therefore, the matrix $H$ has free variables which can be
used in order to satisfy some supplementary requirements for the structure
of the matrix $H$ ensuring that system  \rf{eq1}, \rf{eq2} has exactly $n-r$
zeros.  It is known (Smagina,1990) that the following condition
\be
detHB \neq 0
\label{eq27}
\ee ensures that system  \rf{eq1}, \rf{eq2} possesses the maximal number
of finite zeros.  We transform  system  \rf{eq1}, \rf{eq2} by means of the
above nonsingular state transformation with the matrix  $\bar{N} =NM$. The
transformed output matrix has the form
\be \bar{H} = H\bar{N}^{-1} = \left[
\bar{H}_1, \bar{H}_2 \right]
\label{eq28}
\ee
with $r\times r$  block $\bar{H}_2$.

\begin{ass}\label{A2}
$detHB \neq 0 $ if and only if
\be
det\bar{H}_2 \neq 0.
\label{eq29}
\ee
\end{ass}

{\it Proof.} \rm The following equalities take place
$$ det(HB) = det(H\bar{N}^{-1}\bar{N}B) = det \left( [\bar{H}_1, \bar{H}_2]
\left [ \begin{array}{c} O\\ \tilde{B}_2\end{array} \right ] \right ) =
det(\bar{H}_2\tilde{B}_2)
$$
Since $\tilde{B}_2$  is $r\times r$  nonsingular matrix, then
$det(\bar{H}_2\tilde{B}_2)  \neq 0 $ if and only if $\bar{H}_2$ is a nonsingular
matrix. {\it Q.E.D.}

Restriction \rf{eq29} will be used for the reduction of {\it Problem~\ref{P1}}
to the regular pole assignment problem.
 \begin{th}\label{Th3}
If we assign  $det\bar{H}_2 \neq 0$ then the zero assignment problem in
 \rf{eq1}, \rf{eq2} is reduced to the pole assignment problem in the regular
state-space system of order $n-r$ with $r$ inputs
\be
\dot{\eta } = \bar{A}_{11}\eta -
\bar{A}_{12}\bar{H}_2^{-1}\nu
\label{eq30}
\ee
by using the closed-loop proportional regulator
\be
 \nu = K\eta
\label{eq31}
 \ee where the $(n-r)\times
(n-r)$  matrix $\bar{A}_{11}$ and $(n-r)\times r$ matrix $\bar{A}_{12}$ are
defined from \rf{eq11} with $B_i = \tilde{B}_i, i=1,2$; $r\times (n-r)$  matrix
$K$ coincides with the matrix $\bar{H}_1$.  \end{th}

{\it Proof.} If $det\bar{H}_2 \neq 0$, then  condition \rf{eq27} takes place
in \rf{eq1}, \rf{eq2}  in virtue of {\it Assertion~\ref{A2}}. Therefore, the
polynomial $\psi_a(s)$ in \rf{eq15} has the degree $\mu = n-r$. Using the
determinant expansion formula (Gantmacher, 1959), we can calculate the
determinant of the block matrix in \rf{eq15} with $det{H}_2 \neq 0$
$$
det\left[ \begin{array}{cc} sI_{n-r}-\bar{A}_{11} & -\bar{A}_{12}\\
       \bar{H}_1 & \bar{H}_2  \end{array} \right] = det\bar{H}_2
  det(sI_{n-r} - $$
\be
  - \bar{A}_{11} + \bar{A}_{12}\bar{H}_2^{-1} \bar{H}_1 =
  \psi_a(s)
\label{eq32}
\ee
Since the matrix $\bar{H}_2$ does not depend on $s$, then we can change  pole
assignment condition \rf{eq32} as follows
$$
det(sI_{n-r} - \bar{A}_{11} + \bar{A}_{12}\bar{H}_2^{-1} \bar{H}_1) =
\psi_a(s).
$$
Denoting $\bar{H}_1 =K$, we obtain the regular pole assignment problem.{\it
Q.E.D.}

The necessary and sufficient solvability condition of the regular pole
assignment problem in system  \rf{eq30}  is the controllability of the pair
 \newline $(\bar{A}_{11}, \bar{A}_{12}\bar{H}_2^{-1})$. We can express
this condition in terms of system \rf{eq1},\rf{eq2}.

\begin{ass}\label{A3} If the pair $(A,B)$ is controllable and $rank B = r$,
then the pair $ (\bar{A}_{11}, \bar{A}_{12}\bar{H}_2^{-1})$ is
controllable.
\end{ass}

{\it Proof.} Let the pair $(A,B)$ be controllable, then the rank equality:
$rank(sI-A,B)=n$ takes place for any complex $s$. Using the strict
equivalence transformations, we have series of rank equalities
$$
n = rank\left (sI - A, \; B \right ) = rank\left ( \left [ \begin{array}{cc}
\bar{N} & O \\ O & I  \end{array} \right ] \left[ sI - A, \; B \right ]
\bar{N}^{-1} \right ) = $$
$$  = rank \left [ \begin{array}{ccc} sI-\bar{A}_{11} &
  -\bar{A}_{12} & O \\ -\bar{A}_{21} & sI-\bar{A}_{22} & \tilde{B}_2
   \end{array} \right ] =
rank \left (sI - \bar{A}_{11},- \bar{A}_{12} \right ) + r = $$
$$  = rank \left ( \begin{array}{cc} sI-\bar{A}_{11} &
 -\bar{A}_{12} \end{array} \right) \left [ \begin{array}{cc}
  I & O \\ O & \bar{H}_2^{-1}  \end{array} \right ]  + r =
  rank \left (sI-\bar{A}_{11}, - \bar{A}_{12} \bar{H}_2^{-1}
   \right ) + r
  $$
Hence, $ rank \left ( sI - \bar{A}_{11}, \; - \bar{A}_{12} \bar{H}_2^{-1}
      \right ) = n-r. $ {\it Q.E.D.}
\begin{rem}\label{R2}
{\it Assertion~\ref{A3}} gives the general solvability conditions of
{\it Problem~\ref{P1}} with $\mu = n-r$: the pair $(A,B)$ is controllable and
$rank B =r$. These conditions coincide with the appropriate conditions
of {\it Theorem~\ref{Th2}}.
\end{rem}
\begin{rem}\label{R3}
It should be noted that the nonsingularity of the block $\bar{H}_2$
ensures the rank fullness of the matrix $H$ because
$r=rank[\bar{H}_1,\bar{H}_2] = rankH\bar{N}^{-1} = rankH$.
\end{rem}
Therefore, the present method guarantees the rank fullness of $H$ in
constructed system \rf{eq1}, \rf{eq2} having $n-r$ preassigned zeros.

In conclusion, we write the general zero assignment algorithm.

1. Verify the controllability of $(A,B)$. If $(A,B)$ is not controllable then
{\it Problem~\ref{P1}} has no solution.

2. Assign a desirable polynomial $\psi_a(s)$ of order $ \mu \leq n-r$ with
zeros $s_i$ satisfying {\it Remark~\ref{R1}}.

3. Calculate $M$ and $N$ (formula \rf{eq9}) for $B_i = \tilde{B}_i,i=1,2$.

4. Calculate $\bar{A}_{11}$ and $\bar{A}_{12}$ (formula \rf{eq11}) for $B_i =
\tilde{B}_i,i=1,2$.

5. If $\mu < n-r$, then solve the pole assignment problem in descriptor system
   \rf{eq18} otherwise assign the nonsingular block $\bar{H}_2$ and
    solve the regular pole assignment problem in system \rf{eq30}.

6. Calculate $H$ from \rf{eq16} or \rf{eq23}.

\section{\small{EXAMPLES}}
To illustrate the main results we consider system \rf{eq1} with $n=4$, $r=2$
and
\be  A = \left[ \begin{array}{cccc}
        2  & 1  &  0  & 0 \\
        0  & 1  &  0  & 1 \\
        0  & 2  &  0  & 0 \\
        1  & 1  &  0  & 0 \end{array} \right],\qquad
  B = \left[\begin{array}{cc}   1 & 0\\
                                0 & 0\\
                                0 & 1\\
                                0 & 1     \end{array}\right].
\label{eq33}
\ee

{\it EXAMPLE} 1. It is necessary to find $2\times 4$ output matrix $H$
which guarantees $\mu =2$ finite zeros: $s_1=-1$, $s_2=-2$  in  system
\rf{eq1}, \rf{eq2} with $A$, $B$ from  \rf{eq33}, i.e.
\be
 \psi_a(s) = s^2 + 3s + 2.
\label{eq34}
\ee

We can reduce the zero assignment to the pole assignment problem in a regular
 system because
$\mu$ is equal to the maximal number of zeros in a system.  Pair \rf{eq33} is
controllable, $rankB =2$ and the eigenvalues of $A$ don't coincide with the
preassigned zeros, therefore, {\it Problem~\ref{P1}} has a solution. As the two
last rows of the matrix $B$ are linearly dependent, then we should rearrange
the rows of $B$ by the permutation matrix $$  M= \left[ \begin{array}{cccc} 0  & 0
        &  1  & 0 \\ 0  & 1  &  0  & 0 \\
        1  & 0  &  0  & 0 \\
        0  & 0  &  0  & 1 \end{array} \right],
$$
therefore,
\be     MB = \left[\begin{array}{cc}  0 & 1\\ 0 & 0\\ 1 & 0\\ 0 & 1
          \end{array}\right] ,\;\;
\tilde{B}_1 = \left[\begin{array}{cc}  0 & 1\\ 0 & 0\end{array}\right],\;\;
\tilde{B}_2 = \left[\begin{array}{cc}  1 & 0\\ 0 & 1\end{array}\right].
\label{eq35}
\ee
For calculating $\bar{N}A\bar{N}^{-1}$ we first find
$$ MAM^{T} = \left[ \begin{array}{cccc}
               0  & 2  &  0  & 0 \\
               0  & 1  &  0  & 1 \\
               0  & 1  &  2  & 0 \\
               0  & 1  &  1  & 0 \end{array} \right].
$$
Using the partition of the matrix $MAM^{T} =
\left[ \begin{array}{cc} \bar{A}^*_{11} & \bar{A}^*_{12} \\ \bar{A}^*_{21}&
\bar{A}^*_{22} \end{array} \right]$ with $$ \bar{A}^*_{11} =
\left[\begin{array}{cc}  0 & 2\\ 0 & 1\end{array}\right], \bar{A}^*_{12} =
 \left[\begin{array}{cc}  0 & 0\\ 0 & 1\end{array}\right], \bar{A}^*_{21} =
\left[\begin{array}{cc}  0 & 1\\ 0 & 1\end{array}\right], \bar{A}^*_{22} =
 \left[\begin{array}{cc}  2 & 0\\ 1 & 0\end{array}\right] ,
$$
we calculate the blocks $\bar{A}_{11}, \bar{A}_{12}$ from formulas
\rf{eq11} with $A_{ij} =$ $ \bar{A}^*_{ij}$ and $B_i =\tilde{B}_i$, $i=1,2 $
from \rf{eq35}
$$ \bar{A}_{11}= \bar{A}^*_{11}
-\tilde{B}_1\tilde{B}_2^{-1}\bar{A}^*_{21} = \left[\begin{array}{cc}  0 & 1\\ 0
& 1\end{array}\right], $$
$$
\bar{A}_{12}=\bar{A}^*_{11}\tilde{B}_1\tilde{B}_2^{-1}+\bar{A}^*_{12}
          - \tilde{B}_1\tilde{B}_2^{-1}\bar{A}^*_{22} =
\left[\begin{array}{cc} -1 & 0\\ 0 & 1\end{array}\right]
$$
Assigning the $2\times 2$ block $\bar{H}_2$ in the $2\times 4$ matrix
 $\bar{H} = $ $[\bar{H}_1, \bar{H}_2]$ as unity matrix, we arrive at the pole
assignment problem in the system
\be \dot{\eta } =
\left[\begin{array}{cc}  0 & 1\\ 0 & 1\end{array}\right] \eta -
\left[\begin{array}{cc} -1 & 0\\ 0 & 1\end{array}\right] \nu
 \label{eq36}
 \ee
by the feedback state regulator
 \be \nu = \left[\begin{array}{cc}  k_{11} &
k_{12}\\ k_{21} & k_{22} \end{array}\right]\eta.
 \label{eq37}
 \ee
This problem has a solution because system \rf{eq36} is controllable.
It is evident that the matrix $K$ has free elements, therefore, a family
of  matrices ${\cal K}$ will satisfy a solution. We assign $k_{11} =1$,
$k_{12} =0$  and calculate
$$ K =  \left[\begin{array}{cc}  1 & 0\\ 6 & 5\end{array}\right].  $$
Since $K =\bar{H}_1$ and $\bar{H}_2 =I_2$, then
$$
 \bar{H} = \left[\begin{array}{cccc}  1 & 0& 1 & 0\\ 6 & 5 & 0 & 1
            \end{array}\right].
$$
Using \rf{eq23}, we calculate
$$
H = \bar{H}\bar{N} = \left[ \begin{array}{cccc}  1 & 0 & 1 & -1\\ 0 & 5 & 6 &-5
                            \end{array}\right].
$$
The reliability of this result is verified by calculating  $detP(s)$ for $A$, $B$
from \rf{eq33} and above $H$.

{\it EXAMPLE} 2. It is desirable to design output \rf{eq2} for system
\rf{eq1} with $A$ and $B$ from \rf{eq33} so that constructed system has
the one finite zero $s_1 = -1$. In this case $\mu =1$, $\psi_a(s) = s+1$.
It follows from {\it Theorem~\ref{Th2}} that the problem has a solution because
$(A,B)$ is controllable, $rank B =2$ and $s_1 \neq \lambda_i(A)$.

By substituting the calculated matrices  $\bar{A}_{11}$, $\bar{A}_{12}$
from {\it Example~}1 in \rf{eq17} we construct the descriptor system
\rf{eq18} of order $4$ with two inputs and
$$
E =  \left[\begin{array}{cc}  I_2 & O\\ O & O\end{array}\right],
F =  \left[\begin{array}{cc}  \bar{A}_{11} & \bar{A}_{12} \\ O &
     O\end{array}\right],
 G =  \left[\begin{array}{cc}  O \\ -I_2 \end{array}\right].
$$
It is necessary to shift one finite pole \footnote{A second pole of the
closed-loop system is situated in $\infty $.} of this system to $-1$ by
proportional regulator\rf{eq19} with $ K = \left[\begin{array}{cccc}
  k_{11} & k_{12} & k_{13} & k_{14} \\ k_{21} & k_{22} & k_{23} &
k_{24}\end{array}\right] $. We can calculate some $K$
from the family of matrices ${\cal K}$ satisfying condition \rf{eq20} with
the abovementioned $E$, $F$, $G$
$$  K = \left[\begin{array}{cccc} 0 & 2 & 0 & 1\\
1 & 0 & 0 & 0 \end{array}\right]
$$
Since $K = \bar{H}$ then we find the matrix $H$ using formula
\rf{eq23}
$$
H =  \left[\begin{array}{cccc} 0 & 2 & 0 & 1\\
0 & 0 & 1 & -1 \end{array}\right].
$$ The calculation of $P(s)$ confirms the
correctness of the solution.

\section{\small{CONCLUSION}}
We have presented the original method which reduces the zero assignment problem
to the pole assignment in a singular or regular dynamical systems.
The equivalence between the zero and the pole assignment problems is proved.
The solvability condition of the zero assignment problem in terms of
the pair $(A,B)$ is obtained.


\begin{thebibliography}{99}

\bibitem{Arm:84}
V.A. Armentano.
\newblock Eigenvalue placement for generalized linear systems.
\newblock {\em System and Control letters}, {\bf 4},199--202, 1984.

\bibitem{Bla:90}
F.~Blanchini.
\newblock Controllability analysis and eigenvalue assignment for generalized
  state-space systems.
\newblock {\em System and Control letters}, {\bf 15},285--293, 1990.

\bibitem{Chu:88}
K.W.E. Chu.
\newblock Controllability condensed form and a state feedback pole assignment
  algorithm for descriptor systems.
\newblock {\em IEEE Trans. Autom.Control}, {\bf AC-33},366--370, 1988.

\bibitem{Cob:81}
D.J. Cobb.
\newblock Feedback and pole assignment in descriptor variable systems.
\newblock {\em Int.J.Control}, {\bf 33},1135--1146, 1981.

\bibitem{DavWan:74}
E.J. Davison and S.~Wang.
\newblock Property and calculation of transmission zeros of linear
  multivariable systems.
\newblock {\em Automatica}, {\bf 10},643--658, 1974.

\bibitem{Gan:59}
F.R. Gantmacher.
\newblock {\em Theory of matrices}.
\newblock Chelsea Publishing Co., New York, 1959.

\bibitem{MacKar:76}
A.G.J. MacFarlane and N.~Karcanias.
\newblock Poles and zeros of linear multivariable systems: a survey of the
  algebraic, geometric and complex variable theory.
\newblock {\em Int.J.Control}, {\bf 24},33--74, 1976.

\bibitem{OzcLew:90}
K.~Ozcaldiran and F.L. Lewis.
\newblock Regularizability for singular system.
\newblock {\em IEEE Trans. Autom.Control}, {\bf AC-35},1156--1160, 1990.

\bibitem{PuHaFr:87}
A.C. Pugh, G.E. Hayton, and P.~Fretwell.
\newblock Transformation of matrix pencil and implication in linear system
  theory.
\newblock {\em Int.J.Control}, {\bf 45},529--548, 1987.

\bibitem{Ros:70}
H.H. Rosenbrock.
\newblock {\em State-space and multivariable theory}.
\newblock Nelson, London, 1970.

\bibitem{Sma:84}
Ye.M. Smagina.
\newblock Design of multivariable system with assign zeros.
\newblock {\em Dep.in All-Union Inst. of Sci. and Techn. Inf., Moscow,},
  {\bf No.8309-84},1--15, 1984.

\bibitem{Sma:86}
Ye.M. Smagina.
\newblock Security of specified zeros of multivariable system.
\newblock In A.S.Vostrikov, editor, {\em Automatic control of opbjects with
  varing characteristics}, pages 145--151. NETI, Novosibirsk,Russia, 1986.

\bibitem{Sma:87}
Ye.M. Smagina.
\newblock Computing and specification of zeros in a linear multi-dimension
  systems.
\newblock {\em Avtomatika i telemekhanika}, {\bf 12},165--173, 1987.

\bibitem{Sma:91}
Ye.M. Smagina.
\newblock A method of designing of observable output ensuring given zero
  location.
\newblock {\em Problems of Control and Information Theory}, {\bf 20},299--307,
 1991.

\bibitem{Sma:96}
Ye.M. Smagina.
\newblock The relationship between the transmission zero assignment problem a
  modal control method.
\newblock {\em Journal of Computer and System Science Intern.}, {\bf 35},39--
  1996.

\bibitem{SyrLew:93}
V.L. Syrmos and F.L. Lewis.
\newblock Transmission zero assignment using semistate description.
\newblock {\em IEEE Trans. Autom.Control}, {\bf AC-38},1115--1120, 1993.

\bibitem{VeLeKa:81}
G.C. Verhgese, B.C. Levy, and T.~Kailath.
\newblock A generalised state-space for singular systems.
\newblock {\em IEEE Trans. Autom.Control}, {\bf AC-26},811--831, 1981.

\bibitem{YipSin:81}
E.L. Yip and R.F. Sincovec.
\newblock Solvability, controllability and observability of continuous
  descriptor systems.
\newblock {\em IEEE Trans. Autom.Control}, {\bf AC-26},702--707, 1981.

\end{thebibliography}
\end{document}